\documentclass[letterpaper, 10 pt, conference]{ieeeconf}
\IEEEoverridecommandlockouts
\newcommand{\calv}{\mathcal{V}}
\newcommand{\cale}{\mathcal{E}}
\newcommand{\calb}{\mathcal{B}}
\newcommand{\lhblue}[1]{{\color{black}{#1}}}
\newcommand{\lhred}[1]{{\color{black}{#1}}}
\newcommand{\lhorange}[1]{{\color{black}{#1}}}
\newcommand{\lhgreen}[1]{{\color{black}{#1}}}
\newcommand{\lhviolet}[1]{{\color{black}{#1}}}
\newcommand{\lhteal}[1]{{\color{black}{#1}}}
\newcommand{\lholive}[1]{{\color{black}{#1}}}
\newcommand{\lhmag}[1]{{\color{black}{#1}}}
\newcommand{\icl}[1]{{\color{black}{#1}}}

\usepackage{amsmath}
\usepackage{amsfonts}
\usepackage{algorithmic}
\usepackage{textcomp}
\usepackage{url}
\allowdisplaybreaks

\usepackage{graphicx}
\graphicspath{ {./Images/} }

\usepackage{siunitx}
\DeclareSIUnit\year{yr}
\NewDocumentCommand\SIpi{}{\text{\ensuremath{\pi}}}
\usepackage[dvipsnames]{xcolor}

\usepackage[draft]{changes}

\newtheorem{theorem}{Theorem}

\newtheorem{lemma}[theorem]{Lemma}

\newtheorem{remark}{Remark}

\usepackage[short]{optidef}

\usepackage{tikz}
\usetikzlibrary{positioning}
\usepackage[european,straightvoltages]{circuitikz}

\usepackage{hyperref}
\hypersetup{
	citebordercolor = {.214 .214 .214},
	linkbordercolor = {.214 .214 .214},
	menubordercolor = {.214 .214 .214},
	urlbordercolor = {.214 .214 .214},
	pdfborderstyle={/S/U/W 1},
}
\newcommand{\ic}{\textcolor{black}}

\begin{document}

\title{\LARGE \bf Inverse Optimal Control and Passivity-Based Design for Converter-Based Microgrids}

\author{Liam Hallinan, Jeremy D. Watson, Ioannis Lestas %
\thanks{\lhteal{L. Hallinan and I. Lestas are with with the Department of Engineering, University of Cambridge, Trumpington Street, Cambridge, CB2 1PZ, United Kingdom. Emails:
        {\tt\small <lh706, icl20>@cam.ac.uk}}}%
\thanks{\lhteal{J. Watson is with the Department of Electrical and Computer Engineering, University of Canterbury, Christchurch, Canterbury, New Zealand. Email:
        {\tt\small jeremy.watson@canterbury.ac.nz}}}%
}
\maketitle
\thispagestyle{empty}
\pagestyle{empty}

\begin{abstract}
Passivity-based approaches have been suggested as a solution to the problem of decentralised \lhblue{control} design in many multi-agent network control problems due to the \lholive{plug-and-play functionality they provide}. However, it is not clear if these controllers are optimal at a network level due to their inherently local formulation, with designers often relying on heuristics to achieve desired global performance. On the other hand, solving for an optimal controller is not guaranteed to produce a passive system. In this paper, we address these dual problems by using inverse optimal control theory to formulate a set of sufficient local conditions, which when satisfied ensure that the resulting decentralised control policies are the solution to a network optimal control problem, while at the same time satisfying appropriate passivity properties. These conditions are then reformulated into a set of linear matrix inequalities (LMIs) which can be used to yield such controllers for linear systems. The proposed approach is demonstrated through a DC microgrid case study. The results substantiate the feasibility and efficacy of the presented method.
\end{abstract}

\section{Introduction}
Passivity-based approaches are widely used in the control of dynamic multi-agent networks as they enable the design of decentralised controllers which have desired plug-and-play characteristics, i.e. stability is ensured when subsystems are added/removed from the network. \cite{yao_PassivitybasedControlSynchronization2009, arcak_PassivityDesignTool2007}. Such approaches have been applied to networks involving congestion control \cite{wen_UnifyingPassivityFramework2004} and chemical and biological systems \cite{yao_PassivitybasedControlSynchronization2009}. One particular case of interest is that of power systems. With the rise in distributed energy resources (DER) such as solar photovoltaics and battery storage, the traditional top-down structure of the grid \lhorange{is transitioning} towards a decentralised network of microgrids \lhred{with dynamics dominated by power-electronic converters} \cite{schiffer_survey_2016}. Passivity-based approaches to grid\lhblue{-}forming converter control have been proposed to guarantee decentralised stability \lhred{in both AC \cite{watson_scalable_2021, strehle_unified_2021} and DC microgrids \cite{DecentralizedStabilityConditionslaib2022, PassivitybasedApproachVoltagenahata2020}}. However, the formulation of these controllers typically have no direct relation to the wider network, meaning it is not clear how to measure the performance at a global network level.

In particular, it is desirable in an interconnected system, such as a microgrid, to have a comprehensive metric that assesses the performance of a controller at the network level in order to identify the optimal solution. It was shown in \cite{burger_DualityNetworkTheory2014, sharf_AnalysisSynthesisMIMO2019} that a duality exists between passivity-based cooperative control problems and network optimisation problems where the steady-state equilibria of such systems are solutions to certain optimal network problems, but these results do not consider the problem of finding a control law that minimises a time-domain cost functional. While the passive controllers from the literature discussed above may be optimal at a local level (e.g. the mixed $H_{\infty}$/passive controller from \cite{watson_scalable_2021}), the tuning of these controllers is often based on heuristics and lacks a systematic way of achieving network-wide performance. Here, we \lhblue{propose} an approach \lhred{for addressing} this problem.

The contributions of this paper are as follows:
\begin{itemize}
\item We \lhred{show} that when a set of sufficient conditions are satisfied at each local node in a network, \lhred{then decentralised control policies can be constructed that solve a network wide optimal control problem and at the same time passivate the subsystem dynamics at each node. This leads to a plug-and-play functionality whereby stability is guaranteed as the network is modified while retaining the network wide optimal control formulation.}
\item Next, we provide a methodology for finding such controllers for linear systems \lhred{at each node}, which will involve imposing a set of sufficient constraints on the local node dynamics in the form of a set of linear matrix inequalities (LMIs).
\item Within this framework, it is shown that tuning can be achieved by altering a design matrix for each local controller to account for different design goals, with the new controllers simultaneously satisfying the optimality and passivity requirements outlined above.
\item This approach is tested and verified using a DC microgrid case study.
\end{itemize}


The paper is organised as follows. Section \ref{sec_preliminaries} introduces the network model and reviews both passivity theory and inverse optimal control. These \lhblue{concepts} are combined in Section \ref{sec_main_theory} to design a class of distributed controllers that 
\ic{passiviate the local subsystems}
and also \ic{solve} 
a network-wide optimal control problem. A DC microgrid case study \lhblue{is} analysed in Section \ref{sec_case_studies}, before \lhteal{concluding} in Section \ref{sec_conclusion}.

\section{Preliminaries} \label{sec_preliminaries}

\subsection{Notation and Definitions}
The set $\mathbb{R}_{>0}$ denotes the set of positive real numbers. The matrix-weighted norm of a vector $x \in \mathbb{R}^n$ is given by $\lVert x \rVert_R^2 = x^T R x$ for $R \in \mathbb{R}^{n \times n}$. $I_n$ denotes the n-dimensional identity matrix, with \ic{the} subscript omitted when \ic{the} dimension is obvious. The Kronecker product is denoted by $\otimes$. The direct sum of a set of indexed matrices $B_k$, with $k$ \lhblue{an} element of some ordered index set $\mathcal{A}$ is denoted by $\oplus_{k \in \mathcal{A}} B_k$. The composite vector constructed from a set of indexed vectors $a_k$, where $k$ is an element of some ordered index set $\mathcal{A}$ is denoted $[a_k]_{k \in \mathcal{A}}$. \lholive{A block-diagonal matrix with diagonal entries $X_1, X_2, \ldots, X_n$ is denoted $\mathrm{blockdiag} \begin{bmatrix}
X_1 & X_2 & \ldots & X_n
\end{bmatrix}$.}

\subsection{Network Model} \label{sec_network_model}
Consider a network modelled as a directed graph $\mathcal{G}(\mathcal{V},\mathcal{E})$, with $\mathcal{V} = \{ \nu_1, \nu_2, \ldots, \nu_{|\calv|} \}$ the set representing \ic{the 
nodes} and $\mathcal{E} =  \{ \varepsilon_1, \varepsilon_2, \ldots, \varepsilon_{|\cale|} \}$ is the set representing \ic{the} 
edges \ic{of the graph}. Assigning each edge in $\mathcal{E}$ an arbitrary direction, we denote $N_i^+ \subset \mathcal{E}$ as the set of edges that have node $i$ as their sink and $N_i^- \subset \mathcal{E}$ as the set of edges that have node $i$ as their source. The edge $\varepsilon_k \in \mathcal{E}$ can equivalently be identified as the couple $\varepsilon_k \equiv (\nu_i,\nu_j) \in \mathcal{E} \subset \mathcal{V} \times \mathcal{V}$, indicating that $\varepsilon_k$ connects node $\nu_i \in \calv$ to node $\nu_j \in \mathcal{V}$. Let $\calb \in \mathbb{R}^{|\calv| \times |\cale|}$ denote the incidence matrix of the graph, with $\calb_{jk} = 1$ if node $j \in \mathcal{V}$ is the sink node of edge $k \in \mathcal{E}$, $ \calb_{jk} = -1$ if node $j \in \mathcal{V}$ is the source node of edge $k \in \mathcal{E}$, or $0$ otherwise. Let $\calb_p$ denote the matrix $\calb \otimes I_p$.

\lhred{For each} node $i \in \mathcal{V}$ \lhred{we associate} controllable dynamics that depend non-linearly on local state $x_i \in \mathbb{R}^n$ and affinely on local control variable $u_i \in \mathbb{R}^m$ and \ic{input $w_i \in \mathbb{R}^p$ (defined in \eqref{eq_kcl} in terms of edge variables)} as follows:
\begin{equation} \label{eq_node_sys}
\Sigma_i:
\begin{cases}
\dot{x}_i & = f_i(x_i) + B_{ui}(x_i) u_i + B_i(x_i) w_i \\
y_i & = h_i(x_i, w_i)
\end{cases}
\end{equation}
Here, \lhred{$f_i: \mathbb{R}^n \rightarrow \mathbb{R}^n$, $B_{ui}: \mathbb{R}^n \rightarrow \mathbb{R}^{n \times m}$ and $B_{i}: \mathbb{R}^n \rightarrow \mathbb{R}^{n \times p}$ are Lipschitz continuous. In addition, the system associated with} $i \in \mathcal{V}$ produces output $y_i \in \mathbb{R}^p$ \lhred{and $h_i: \mathbb{R}^n \times \mathbb{R}^p \rightarrow \mathbb{R}^p$ is a continuous vector-valued function}. \lhteal{We} \ic{focus on} 
the case of \textit{decentralised} control, meaning that \ic{the control} variable $u_i$ depends on local measurements only, with no requirement for inter-node communication. \lhblue{In the paper, we will consider static state-feedback policies of the form} \lhred{$u_i = k_i(x_i)$, where $k_i: \mathbb{R}^n \rightarrow \mathbb{R}^m$ is continuous}.

Similarly, \lhred{with} each edge $j \in \mathcal{E}$ \lhred{we associate} uncontrollable dynamics \lhred{that} depend non-linearly on \lhred{the} local state $x_j \in \mathbb{R}^r$ and affinely on the \ic{input $w_j\in \mathbb{R}^p$ (defined in \eqref{eq_kvl} in terms of node variables)} as follows:
\begin{equation} \label{eq_edge_sys}
\Sigma_j:
\begin{cases}
\dot{x}_j & = f_j(x_j)  + B_j(x_j) w_j \\
y_j & = h_j(x_j,w_j)
\end{cases}
\end{equation}
Here, \lhred{$f_j: \mathbb{R}^r \rightarrow \mathbb{R}^r$ and $B_{j}: \mathbb{R}^r \rightarrow \mathbb{R}^{r \times p}$ are Lipschitz continuous. In addition, the system associated with } $j \in \mathcal{E}$ produces output $y_j \in \mathbb{R}^p$ \lhred{and $h_j: \mathbb{R}^r \times \mathbb{R}^p \rightarrow \mathbb{R}^p$ is continuous.}

Motivated by electrical circuit theory, the interconnection between the node and edge dynamics is characterised by the following relations (which respectively correspond to Kirchhoff's Current Law and Kirchhoff's Voltage Law):
\begin{subequations} \label{eq_inout_interconnection}
\begin{align}
& w_i  =  \sum_{(l,m) \in N_i^-} y_{lm} - \sum_{(l,m) \in N_i^+} y_{lm}, & & i \in \mathcal{V} & \label{eq_kcl} \\
& \qquad \qquad \quad w_{lm} = - y_{l} +  y_{m}, & & (l,m) \in \mathcal{E}   \label{eq_kvl}
\end{align}
\end{subequations}
\lhred{Here, $w_{lm}$ and $y_{lm}$ \lhteal{respectively} denote the input and output of system \eqref{eq_edge_sys} \lhorange{associated with} the edge $(l,m) \in \cale$.}

Now, using the above relations, two composite subsystems can be created representing the node and edge dynamics. Letting $x_{\calv} = [x_i]_{i \in \mathcal{V}}$, $f_{\calv}(x_{\calv}) = [f_i(x_i)]_{i \in \mathcal{V}}$, $w_{\calv} = [w_i]_{i \in \mathcal{V}}$, $\hat{u} = [u_i]_{i \in \mathcal{V}}$, $B_{\calv}(x_{\calv}) = \oplus_{i \in \mathcal{V}} B_i(x_i)$ and $B_u(x_{\calv}) = \oplus_{i \in \mathcal{V}} B_{ui}(x_i)$, we have the composite \lhorange{block diagonal} node dynamics
\begin{align}
& \Sigma_{\calv}: &  &\dot{x}_{\calv} = f_{\calv}(x_{\calv}) + B_{\calv}(x_{\calv})w_{\calv} + B_u(x_{\calv})\hat{u} &
\end{align}
Similarly, the composite \lhorange{block diagonal} edge dynamics can be described using $x_{\cale} = [x_j]_{j \in \mathcal{E}}, f_{\cale}(x_{\cale}) = [f_j(x_j)]_{j \in \mathcal{E}}, w_{\cale} = [w_j]_{j \in \mathcal{E}}$ and $B_{\cale}(x_{\cale}) = \oplus_{j \in \mathcal{E}} B_j(x_j)$, giving
\begin{align}
& \Sigma_{\cale}: & & \quad \dot{x}_{\cale} = f_{\cale}(x_{\cale}) + B_{\cale}(x_{\cale})w_{\cale} &
\end{align}
The relations \eqref{eq_inout_interconnection} then reduce to $w_{\calv} =  -\calb_p y_{\cale}$ and $w_{\cale} = \calb_p^T  y_{\calv}$, where $y_{\calv} = [y_i]_{i \in \mathcal{V}}$ and $y_{\cale} = [y_j]_{j \in \mathcal{E}}$. Therefore, the whole network can be described by $\dot{\hat{x}} = \hat{f}(\hat{x}) + \hat{B}(\hat{x}) \hat{u}$, \lhred{with $\hat{x} = [x_{\calv}^T \ x_{\cale}^T]^T$ and $\hat{f}(\hat{x})$ and $\hat{B}(\hat{x})$ given by the following expression:}
\begin{equation} \label{eq_full_network}
\begin{bmatrix} \dot{x}_{\calv} \\ \dot{x}_{\cale} \end{bmatrix} =
\begin{bmatrix} f_{\calv}(x_{\calv}) - B_{\calv}(x_{\calv}) \calb_p y_{\cale} \\ f_{\cale}(x_{\cale}) + B_{\cale}(x_{\cale}) \calb_p^T  y_{\calv} \end{bmatrix} +
\begin{bmatrix}
B_u(x_{\calv}) \\ 0
\end{bmatrix} \hat{u}
\end{equation}

\begin{figure}[t]
\centering

\begin{tikzpicture}

\node[draw, circle,
			minimum size=0.6cm]
			(neg) at (0,0) {};
			\node at (neg.center){ $-1$};
\node[draw, rectangle,
			right=of neg,
			minimum width=0.75cm, 	minimum height=0.75cm]
			(sys_1) {$\Sigma_{\calv}$};
\node[draw, rectangle,
			below=of sys_1,
			minimum width=0.75cm, 	minimum height=0.75cm]
			(sys_2) {$\Sigma_{\cale}$};
\node[draw, circle,
			minimum size=0.6cm]
			(B) at (0,-1) {};
			\node at (B.center){ $\calb_p$};
\node[draw, circle,
			right=of sys_1,
			yshift= -1cm,
			minimum size=0.6cm]
			(Bt)   {};
			\node at (Bt.center){ $\calb_p^T$};
			
\draw[-stealth] (neg.east) -- (sys_1.west)
    node[midway,above]{$w_{\calv}$};
\draw[-stealth] (sys_1.east) -| (Bt.north)
    node[pos=0.25,above]{$y_{\calv}$};
\draw[-stealth] (Bt.south) |- (sys_2.east)
    node[pos=0.75,above]{$w_{\cale}$};
\draw[-stealth] (sys_2.west) -| (B.south)
    node[pos=0.25,above]{$y_{\cale}$};
\draw[-stealth] (B.north) -- (neg.south);

\end{tikzpicture}

\caption{Negative feedback interconnection of the node system $\Sigma_{\calv}$ with state $x_{\calv}$ and the edge system $\Sigma_{\cale}$ with state $x_{\cale}$.}
\label{fig_interconnection}
\end{figure}
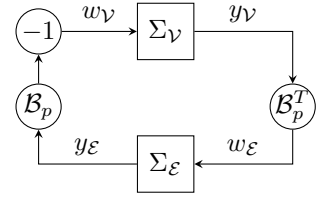

We therefore see that the system \eqref{eq_full_network} in closed-loop is a negative feedback interconnection of two aggregate systems as demonstrated in Figure \ref{fig_interconnection}: the node subsystem $\Sigma_{\calv}$ with state $x_{\calv}$ and the \ic{edge subsystem $\Sigma_{\cale}$} with state $x_{\cale}$.

\subsection{Review of Passivity Theory}
As \lhblue{described} in \cite{yao_PassivitybasedControlSynchronization2009, watson_scalable_2021}, a local state feedback controller can be designed to achieve certain passivity properties \cite[Chapter 6]{khalil_nonlinear_2013} of the closed-loop node subsystem. A similar controller will be designed here.

Strict passivity \cite[Definition 6.3]{khalil_nonlinear_2013} of $i \in \mathcal{D}$ is satisfied when there exists a positive definite \textit{storage function} $V_i(x_i)$ such that
\begin{equation} \label{eq_strict_passivity}
w_i^Ty_i \leq \dot{V}_i(x_i) + \psi_i(x_i)
\end{equation}
where $\psi_i(x_i)$ is a positive definite function. For linear systems $i \in \calv$ of the form
\begin{equation} \label{eq_linearised_sys}
\begin{aligned}
\dot{x}_i &= A_i x_i + B_{ui} u_i + B_i w_i \\
y_i & = C_i x_i
\end{aligned}
\end{equation}
\lhred{with linear controllers $u_i = K_i x_i$, $K_i \in \mathbb{R}^{m \times n}$} we can take $V_i(x_i) = \frac{1}{2}x_i^TP_ix_i$ and $\psi_i(x_i) = \frac{1}{2}x_i^T \Gamma_i^{-1} x_i$ for some matrices satisfying $P_i = P_i^T > 0$ and $\Gamma_i = \Gamma_i^T> 0$. Then, \eqref{eq_strict_passivity} can equivalently be written as the following bilinear matrix inequality (BMI) in $P_i$ and $K_i$ \cite{watson_scalable_2021,khalil_nonlinear_2013}:
\begin{equation}
\label{eq_solve_LMI}
\begin{bmatrix}
\tilde{A}^T P_i  + P_i \tilde{A} + \Gamma_i^{-1} & P_i B_{i}- C_i^T \\
B_{i}^T P_i - C_i & 0
\end{bmatrix} \leq 0
\end{equation}
where $\tilde{A} = A_i + B_{ui}K_i$. This BMI can then be converted to an LMI using simple transformations (as \lhred{in} the proof of Lemma \ref{lemma_lmi} in Appendix \ref{app_proof_l3}).

Now, assuming \lhred{the system associated with} each edge $j \in \mathcal{E}$ is \lhviolet{strictly} passive from $w_j$ to $y_j$ with storage function $V_j(x_j)$, and \lhred{the system associated with} each node $i \in \mathcal{V}$ can be rendered \lhviolet{strictly} passive from $w_i$ to  $y_i$ with storage function $V_i(x_i)$, we can then construct a storage function for the entire network $\hat{V}(\hat{x}) =  \sum_{i \in \mathcal{V}} V_i(x_i) + \sum_{j \in \mathcal{E}} V_j(x_j)$. If all systems are linear, this can be written in matrix form as $\hat{V}(\hat{x}) = \frac{1}{2} \hat{x}^T \hat{P} \hat{x}$, or
\begin{equation} \label{eq_vhat_linear}
\hat{V}(\hat{x})  = \frac{1}{2} \begin{bmatrix}
x_{\lhmag{\calv}}^T  \\ x_{\lhmag{\cale}}^T
\end{bmatrix}^T
\begin{bmatrix}
\bigoplus_{i \in \mathcal{V}} P_i & 0  \\
0 &  \bigoplus_{j \in \mathcal{E}} P_{j}
\end{bmatrix}
\begin{bmatrix}
x_{\lhmag{\calv}}  \\ x_{\lhmag{\cale}}
\end{bmatrix}
\end{equation}

The function $\hat{V}(\hat{x})$ can then be used as a Lyapunov function to determine the stability properties of the network dynamics \eqref{eq_full_network}. In addition, \ic{the 
controller $K_i$ can be designed such that} \lhred{the system associated with} each node $i \in \mathcal{V}$ remains \lhviolet{strictly} passive, thus \icl{ensuring} 
stability under changes in network topology \cite{watson_scalable_2021}. However, controllers designed using passivity-based \ic{methods 
do not} necessarily guarantee desirable or optimal performance due to their inherent conservatism. We show here that inverse optimal control theory can be used to address this problem, \ic{i.e. via appropriate conditions on the decentralized controllers 
a network-wide performance metric can be constructed and at the same time the controllers 
\lholive{lead to} \ic{local passivity~properties that facilitate a plug-and-play operation.}}

\subsection{Review of Inverse Optimal Control}
\lhblue{The aim of inverse optimal control is to find the performance metric for which a given controller is the optimal solution}. This can be useful to gain a deeper understanding of the system constraints and how tuning the controller may impact performance.

The problem of reconstructing the objective function from a given controller has a long history, and various formulations have been presented in the literature (see \cite{kalman_WhenLinearControl1964, InverseProblemLinearjameson1973, freeman_InverseOptimalityRobust1996} for some earlier results and \cite{zhang_distributed_2015, jouini_InverseOptimalControl2021} for some more recent accounts). The approach taken in \cite{jouini_cost_2022} for nonlinear systems is reproduced here:
\begin{theorem} \label{theorem_inverse_opt}
Consider the nominal optimal control problem
\begin{mini}
	{u}{\int_0^\infty q(x(s), R) + \lVert u(x(s)) \rVert^2_R ds}{\label{opt_nonlinear}}{}
	\addConstraint {\lholive{ \dot{x}=f(x) + G^T(x)u, \ x(0)=x_0}}
\end{mini}
Here $x \in \mathbb{R}^n$ denotes the state vector, $x(0) = x_0$ is the initial condition and $f(x): \mathbb{R}^n \rightarrow \mathbb{R}^n$ is a continuous nonlinear vector field with $f(0) = 0$. The input matrix $G(x) = \left[ g_1^T(x), \ldots, g_m^T(x) \right]^T \in \mathbb{R}^{m \times n}$ is given by the continuous nonlinear functions $g_i(x): \mathbb{R}^n \rightarrow \mathbb{R}^m, i = 1, \ldots, m$. We have $R = R^T$ as a design matrix, and the function $q(x(s)): \mathbb{R}^n \rightarrow \mathbb{R}_{>0}$ has the condition $q(0) = 0$.

Now, let $V : \mathbb{R}^n \rightarrow \mathbb{R}_{>0}$ a continuously differentiable function associated with a stabilising feedback control law
\begin{equation}
\label{eq_control_law}
u^{\star} = - \tfrac{1}{2} R^{-1} G(x) \nabla V
\end{equation}
where
\begin{equation}
\label{eq_condition}
\nabla V^T \left( f\left(x\right) + G^T\left(x\right)u^{\star}\left(x,R\right) \right) < - \lVert u^{\star} \left(x,R\right)\rVert^2_R.
\end{equation}
Define
\begin{multline}
\label{eq_state_cost_general}
q(x,R) = - \nabla V^T \left(f\left(x\right) + G^T\left(x\right) u^{\star}\left(x,R\right) \right) \\ - \lVert u^{\star} \left( x,R \right) \rVert^2_R.
\end{multline}
Then the following statements hold:
\begin{enumerate}
\item The unique optimal control law is given by \eqref{eq_control_law}.
\item The optimal control problem \eqref{opt_nonlinear} has the optimal value $V(x_0)$.
\end{enumerate}
\end{theorem}

Therefore, \eqref{opt_nonlinear} serves as a performance metric for the system. In addition, as \lhblue{described} in \cite{jouini_cost_2022}, $V(x)$ can be used as a Lyapunov function to prove asymptotic stability\lhviolet{, though additional system structure needs to be exploited when the right-hand-side of \eqref{eq_condition} is not positive definite in $x$}.

For linear systems, the problem reduces to a linear quadratic regulator (LQR) problem \cite{bertsekas_dynamic_2000} of the form:
\begin{mini}
	{u}{\int_0^\infty \lVert x \rVert^2_{Q(R)} + \lVert u \rVert^2_{R} dt}{\label{opt_lqr}}{}
	\addConstraint {\lholive{ \dot{x}=Ax + Bu, \ x(0)=x_0}}
\end{mini}
with state cost matrix $Q(R)= Q(R)^T \geq 0$, input cost matrix $R = R^T > 0$ and initial condition $x_0$. Applying Theorem \ref{theorem_inverse_opt} requires the optimal controller to have an associated control Lyapunov function $V(x) = \frac{1}{2}x^TPx$ with $P = P^T > 0$ and be of the form
\begin{equation} \label{eq_control_law_linear}
u^{\star}(x,R) = - \tfrac{1}{2}R^{-1}B^TPx
\end{equation}
Equation \eqref{eq_state_cost_general} gives an expression for the state cost matrix $Q(R)$ and \eqref{eq_condition} gives the condition $Q(R) > 0$, which together reduce to
\begin{equation}
\label{eq_Q_condition}
Q(R) = \tfrac{1}{4}PBR^{-1}B^TP - \tfrac{1}{2}\left( A^T P + PA\right) > 0.
\end{equation}


\begin{remark} \label{rem_inv_opt_conditions}
Note that Theorem \ref{theorem_inverse_opt} requires the control law to be of the form \eqref{eq_control_law} in order to be the solution to an optimisation problem of the form \eqref{opt_nonlinear}, i.e. given a state feedback control law $u(x)$, we require $R$ and $V(x)$ to exist such that $u(x)$ can be written in this form. This may not be the case for an arbitrary stabilising $u(x)$, \lhblue{meaning an inverse optimal \lhred{control} interpretation may not always be \lhorange{possible}}. For linear systems, a set of necessary and sufficient conditions are \lhblue{specified} in \cite{InverseProblemLinearjameson1973} that must be satisfied for this to be the case.
\end{remark}


\lhmag{For networks \icl{of interconnected systems}, using} a decentralised controller at each node does not necessarily lead to an inverse optimal \lhred{control} interpretation. \lhred{At the same time}, solving a network-wide \lhteal{optimal control problem} is not guaranteed to produce a decentralised set of controllers, nor is it guaranteed to produce \lhred{controllers that lead to passivity properties for the subsystems, and thus a plug-and-play operation.} In the sequel, we will discuss under what conditions a decentralised controller can solve both these problems simultaneously.

\section{Global Optimality from Decentralised Passive Controllers} \label{sec_main_theory}

\lhblue{We see that a passivity based control design leads to a Lyapuonv function for \lhred{the entire} interconnection. At the same time inverse optimal control theory relies on the existence of a Lyapuonv function which satisfies certain additional conditions. We will see in this section that by combining these two approaches} we can derive a set of sufficient local conditions under which a set of decentralised controllers simultaneously yield the optimal solution to a network-wide \lhred{optimal control problem}, while also \lholive{facilitating} plug-and-play \lhblue{operation} \lhteal{when the network is modified}. \lholive{Proofs to \icl{the} results in this section are provided in the Appendix.}

\begin{theorem} \label{theorem_combined}
Consider the negative feedback interconnection network model \lhorange{\eqref{eq_full_network}}. Assume \lhteal{that for each} edge $j \in \mathcal{E}$ \lhteal{the system with input $w_j$ and output $y_j$ is strictly passive} with storage function $V_j(x_j)$. For each node $i \in \mathcal{V}$, assume \lhorange{dynamics \eqref{eq_node_sys} with}
\begin{equation} \label{eq_control_law_local}
u_i = -\tfrac{1}{2}R_i^{-1}B_{ui}^T(x_i) \nabla V_i(x_i)
\end{equation}
with some matrix $R_i = R_i^T > 0$ \lhteal{and} some positive definite function $V_i(x_i)$. Assume that using \eqref{eq_control_law_local} in \eqref{eq_node_sys} generates a local closed-loop system such that the node dynamics are strictly passive from $w_i$ to $y_i$ with storage function $V_i(x_i)$, i.e.
\begin{equation} \label{eq_pass_condition_local}
w_i^Ty_i \geq \dot{V}_i(x_i) + \psi_i(x_i)
\end{equation}
where $\psi_i(x_i)$ is some positive definite function.

If in addition at each node $i \in \mathcal{V}$, the following local condition is satisfied
\begin{multline} \label{eq_condition_local}
- \nabla^T V_i(x_i) f_i(x_i) + \\ \tfrac{1}{4} \nabla^T V_i(x_i) B_{ui}(x_i)R_i^{-1}B_{ui}^T \nabla V_i(x_i) > 0
\end{multline}
then $
\hat{u}^{\star} = \left[ u_i \right]_{i \in \mathcal{V}}
$
is the globally optimal \lhteal{control input} with respect to the optimisation problem \eqref{opt_nonlinear} for system \eqref{eq_full_network} with control cost matrix given by $
R = \oplus_{i \in \mathcal{V}} R_i
$
and state cost function
\begin{multline} \label{eq_state_cost_network}
q(\hat{x}) = \sum_{i \in \mathcal{V}} \left[ - \nabla^T V_i(x_i) f_i(x_i) \right. \\
\left. + \tfrac{1}{4} \nabla^T V_i(x_i) B_{ui}(x_i)R_i^{-1}B_{ui}^T \nabla V_i(x_i) \right] \\
+ \sum_{j \in \mathcal{E}} \left[ - \nabla^T V_j(x_j) f_j(x_j) \right]
\end{multline}

In addition, using \eqref{eq_control_law_local} at each node ensures asymptotic stability of the \lhblue{interconnected system} independent of \lhteal{its} topology.
\end{theorem}


\lholive{
\begin{remark}
It should be noted that modifying the network may lead to a change in the equilibrium point. For linear systems, if the passivity property is satisfied about an equilibrium point, then it is satisfied at any equilibrium point \cite{khalil_nonlinear_2013}\lhmag{. In this case, Theorem \ref{theorem_combined} guarantees stability after network modifications, thus providing plug-and-play functionality.} For nonlinear systems, satisfying the passivity property at any equilibrium point is linked to incremental passivity \cite{stan_AnalysisInterconnectedOscillators2007}.
In addition, the approach outlined here can also be useful in practical designs where a change in \icl{an equilibrium point} does not significantly change the linearisation of the system \icl{at that point}, thus leading to plug-and-play functionality when operating close to this equilibrium point.
\end{remark}
}

It should be noted here (as demonstrated in the proof) that as a result of the local passivity property of each \lhorange{node and edge}, the state cost function \eqref{eq_state_cost_network} is the sum of terms involving local states only. This means that the cost function scales with the size of the network in question and can easily be adjusted to account for new nodes.

Next, a method to find such controllers for linear systems is presented. Here, the model of Section \ref{sec_network_model} is linearised so that it can be written in the form \eqref{eq_linearised_sys}. Similar relationships hold for \eqref{eq_edge_sys}. We let $A_{\calv} = \bigoplus_{i \in \mathcal{V}} A_i$, with similar definitions holding for $B_{\calv}$, $B_{u \calv}$, $C_{\calv}$, $A_{\cale}$, $B_{\cale}$ and $C_{\cale}$. This gives rise to a negative feedback interconnection system between the composite node and edge subsystems  (c.f. \eqref{eq_full_network}) described by:
$
\dot{\hat{x}} = \hat{A}\hat{x} + \hat{B}\hat{u}
$
\lhred{with matrices $\hat{A}$ and $\hat{B}$ given in the following expression:}
\begin{equation} \label{eq_network_linear_full}
\begin{bmatrix} \dot{x}_{\calv} \\ \dot{x}_{\cale} \end{bmatrix} =
\begin{bmatrix}
A_{\calv}  &  - B_{\calv} \calb_p  C_{\cale} \\
B_{\cale} \calb_p^T C_{\calv} &  A_{\calv}
\end{bmatrix} \begin{bmatrix} x_{\calv} \\ x_{\cale} \end{bmatrix} +
\begin{bmatrix}
B_{u\calv}\\
0
\end{bmatrix} \hat{u}
\end{equation}

\begin{lemma} \label{lemma_lmi}
\lhteal{Consider system \eqref{eq_full_network} with linear node and edge dynamics as in \eqref{eq_network_linear_full}. For each $i \in \calv$ with dynamics \eqref{eq_node_sys}, consider matrices $Y_i \in \mathbb{R}^{n \times n}$, $S_i \in \mathbb{R}^{m \times m}$ and $\Gamma_{\lholive{i}} \in \mathbb{R}^{n \times n}$ that satisfy the following set of linear matrix inequalities (LMIs):}
\begin{subequations} \label{opt_lmi}
\begin{gather}
\begin{bmatrix}
Y_iA_i^T + A_iY_i + 2B_{ui}S_iB_{ui}^T  & Y_i & B_{i}- Y_iC_i^T \\
Y_i & - \Gamma_i & 0 \\
B_{i}^T - C_i Y_i& 0 & 0 \end{bmatrix} \leq 0 \label{const_passive} \\
-\tfrac{1}{2} B_{ui} S_i B_{ui}^T - \tfrac{1}{2} \left( A_i Y_i + Y_i A_i^T \right) > 0 \label{const_q_pos_def} \\
Y_i = Y_i^T> 0, \quad \Gamma_i > 0, \quad S_i = S_i^T > 0 \label{const_technicals}
\end{gather}
\end{subequations}
\lhteal{Considering} $P_i = Y_i^{-1}$ and $R_i = -\frac{1}{2}S_i^{-1}$ yields the following distributed local controller for each $i \in \mathcal{V}$:
\begin{equation} \label{eq_uistar}
u_i^{\star} = -\tfrac{1}{2} R_i^{-1}B_{ui}^TP_ix_i
\end{equation}
Then, together with the assumption that the edge subsystems \eqref{eq_edge_sys} for $j \in \cale$, are \lhviolet{strictly} passive from input $w_j$ to output $y_j$, we have that the local controller \eqref{eq_uistar} satisfies the conditions \lhblue{specified} in Theorem \ref{theorem_combined} and therefore:
\begin{enumerate}
\item It generates a strictly passive local closed-loop node subsystem \eqref{eq_node_sys} for $i \in \mathcal{V}$, with input $w_i$ and output $y_i$.
\item The global network-wide controller
\begin{equation} \label{eq_ustar}
\hat{u}^{\star} = \oplus_{i \in \mathcal{V}} \left( -\tfrac{1}{2} R_i^{-1} B_{ui}^T P_i \right) x_{\lhmag{\calv}}
\end{equation}
composed of the distributed controllers \eqref{eq_uistar} is the optimal solution to \eqref{opt_lqr} for system \eqref{eq_network_linear_full} with cost matrices
\begin{subequations} \label{eq_RandQ}
\begin{gather}
R = \oplus_{i \in \mathcal{V}} R_i \\
Q(R) = \lholive{ \mathrm{blockdiag} \begin{bmatrix}
Q_{\calv}(R) & Q_{\cale}
\end{bmatrix} } \label{eq_Q_matrix}
\end{gather}
\end{subequations}
where
\begin{equation}
\begin{aligned} \label{eq_q1q2q3}
Q_{\calv}(R)=& \oplus_{i \in \mathcal{V}} \left(  \tfrac{1}{4} P_i B_{ui}R_i^{-1} B_{ui}^T P_i \right. \\
 & \qquad \left. - \tfrac{1}{2} \left( P_i A_i + A_i^T P_i \right) \right) \\
Q_{\cale} =& -\oplus_{j \in \mathcal{E}} \tfrac{1}{2} \left(P_{jk} A_{jk} + A_{jk}^T P_{jk} \right)
\end{aligned}
\end{equation}
\item The system \eqref{eq_network_linear_full} is asymptotically stable independent of network topology.
\end{enumerate}
\end{lemma}

\begin{remark}
Solving \eqref{const_passive} produces a local controller that generates a strictly passive closed-loop node subsystem, while \eqref{const_q_pos_def} (although it only applies locally) ensures that condition \eqref{eq_condition} is satisfied for the global system. Constraints \eqref{const_technicals} are required so that $P_i$ and $R_i$ are symmetric positive definite, and $\Gamma_i$ is positive definite. Therefore, satisfying all these constraints simultaneously produces a controller that fulfils the passivity property and is also a component of the global optimal network controller $\hat{u}^{\star}$ given by \eqref{eq_ustar} which minimises \eqref{opt_lqr} with $R$ and $Q(R)$ as in \eqref{eq_RandQ}.
\end{remark}

\begin{remark} \label{rem_additional_const}
If desired, additional LMI constraints can be added to \eqref{opt_lmi} to impose desired dynamic behaviour on the system. \lhorange{For example, the} constraint
\begin{equation} \label{eq_lambda}
A_i^TP_i + P_iA_i  + 2P_iB_{ui}S_iB_{ui}^T P_i< \lambda_{\lhmag{i}} P_i
\end{equation}
(converted to an LMI by pre- and post-multiplication by $Y_i$) has the effect of imposing a minimum decay rate of $\lambda_{\lhmag{i}}$ on the storage function $V_i(x_i) = \frac{1}{2}x_i^T P_i x_i$ of the local node subsystem, which can \lhorange{improve} performance by decreasing settling time \lhblue{for the local system} \cite{boyd_linear_nodate}.
\lhblue{As the network Lyapunov function $\hat{V}(\hat{x})$ consists of a sum of the local storage functions (see \eqref{eq_vhat_linear}), this has a direct impact on global performance in cases where the node dynamics dominate and there are weak node-edge interactions.}
\end{remark}

The following result shows that solving for a controller using the method \lhblue{detailed} in Lemma \ref{lemma_lmi} gives rise to the additional benefit that the local controller \eqref{eq_uistar} can easily be tuned using the matrix $R_i$.

Once suitable $P_i$ (and a corresponding $\hat{P}$) have been found for each $i \in \mathcal{V}$, application of Theorem \ref{theorem_inverse_opt} means that new controllers $\bar{u}_i$ can be found via $\bar{u} =  -\frac{1}{2}\bar{R}^{-1}  \hat{B}^T \hat{P} x_{\lhmag{\calv}}$ for any $\bar{R} = \bar{R}^T > 0$. For a decentralised controller, we require a block diagonal $\bar{R}$, so taking $\bar{R} = \oplus_{i \in \mathcal{V}} \bar{R}_i$ with all $\bar{R}_i = \bar{R}_i^T > 0$, gives the local node controller
\begin{equation}
\label{eq_new_K}
\bar{u}_i =  -\tfrac{1}{2}\bar{R}_i^{-1} B_{ui}^T P_i x_i
\end{equation}
and $\bar{u}$ satisfies the global optimality condition with control cost $\bar{R}$, provided $Q\left(\bar{R}\right) \geq 0$.

Thus, the local control cost matrices $\bar{R}_i$ can be used to tune the control matrix gains, while also giving rise to an optimal controller that solves \eqref{opt_lqr}. The following lemma \lhblue{specifies} a condition that ensures the resulting closed-loop node subsystem using controller $\bar{u}_i$ is \lhorange{still strictly passive for new values of $\bar{R}_i$}.

\begin{lemma} \label{lemma_tuning}
Let $R_i$ and $P_i$ denote particular solutions found by solving \eqref{opt_lmi}. Let $\bar{R}_i \in \mathbb{R}^{m \times m}, \bar{R}_i = \bar{R}_i^T > 0$ be a new tuning matrix that gives rise to local controller $\bar{u}_i$ via \eqref{eq_new_K} and global controller $\bar{u}$ via $\bar{u} = \left[ \bar{u}_i \right]_{i \in \mathcal{V}}$. Then, if $\bar{R}_i \leq R_i$,
\begin{enumerate}
\item The local closed-loop node subsystem \eqref{eq_node_sys} with $i \in \mathcal{V}$ is strictly passive from input $w_i$ to output $y_i$.
\item $\bar{u}$ is the optimal solution to the global network LQR problem \eqref{opt_lqr} with cost matrices $\bar{R} = \oplus_{i \in \mathcal{V}} \bar{R}_i$ and $Q(\bar{R})$ given by \eqref{eq_Q_matrix}.
\end{enumerate}
\end{lemma}

\begin{remark} \label{rem_opt_problem}
As a result of Lemma \ref{lemma_tuning}, it may be useful to find a maximal $R_i$ that achieves both the passivity and optimality requirements when solving \eqref{opt_lmi}. \lhblue{With this $R_i$, \lhorange{the control cost can then be lowered by choosing $\bar{R}_i \lhmag{<} R_i$. We also have that the set of such $\bar{u}_i$} form a class of \lhred{controllers that lead to the passivity property for each subsystem which also have an inverse optimal control} interpretation. Using} $S_i = -\frac{1}{2} R_i^{-1}$ in \eqref{opt_lmi}, this can be achieved by letting $S_i = -s_iI_m$ with $s_i \in \mathbb{R}, s_i > 0$ and minimising over the scalar $s_i$, turning \eqref{opt_lmi} into an optimisation problem. Therefore, we have that $\bar{R}_i > R_i$ cannot satisfy both the passivity and optimality conditions simultaneously when $R_i$ is of the form considered here.
\end{remark}

\section{DC Microgrid Case Study} \label{sec_case_studies}
\lhorange{We illustrate our results using a DC microgrid example. This} was chosen as passivity-based approaches \lhteal{in DC microgrids} have been \lhblue{frequently} suggested in the literature \lhorange{where \lhteal{they are} used to provide} a decentralised means for control design without requiring knowledge of the microgrid structure \lhgreen{\cite{ReviewControlDCmeng2017}}. However, the \lhorange{tuning of} such control schemes is often based on heuristics \lhteal{due to the lack of a network-wide performance metric}.

In such a microgrid, we will assume that every bus contains a DC source distributed generation unit (DGU) (such as a solar PV source or a battery), or a load. The DGU units interface with the grid via a DC-DC converter to regulate the DC voltage \cite{DecentralizedStabilityConditionslaib2022, PassivitybasedApproachVoltagenahata2020}.

We will assume that the node set $\mathcal{V}$ will be split between a collection of buses containing DGUs indexed by the set $\mathcal{D}$ and a collection of load buses indexed by the set $\mathcal{L}$, such that $\mathcal{V} =  \{\mathcal{D}, \mathcal{L}\}$.

\begin{figure}[t]
\centering
\begin{circuitikz} \draw
(0,2) to[cute inductor,*-,l=$l_i$] ++(1.5,0)
  to[short,i_=$i_i$] ++(0.25,0)
  to[american resistor,l=$r_i$] ++(1.5,0)
  to[short] ++(0.75,0)
  to[short] ++(1,0)
  to[short,-*,i_=$i_i^0$] ++(2,0)
  to[open,l=$+$] ++(0,-0.25)
  to[open,v^<=$v_i$] ++(0,-1.5)
  to[open,-*,l=$-$] ++(0,-0.25)
  to[short,-*] (0,0)
  to[open,l=$-$] ++(0,0.25)
  to[open,v>=$u_i$] ++(0,1.5)
  to[open,l=$+$] (0,2)
  (3.75,2) to[american resistor,l_=$g_i$] ++(0,-2) node[ground] {}
  (4.75,2) to[capacitor,l=$c_i$] ++(0,-2)
;

\end{circuitikz}
\caption{\lhgreen{Circuit diagram of the buck converter filter.}}
\label{fig_circuit_dc}
\end{figure}
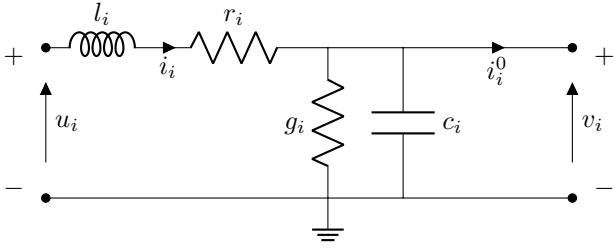

Here, all DGU buses will be modelled with a buck converter producing a \lhorange{controllable} DC voltage. At  the outlet of each converter is \lhblue{an} RLC filter which is connected to the lines of the microgird, as demonstrated in Figure \ref{fig_circuit_dc}. Using Kirchhoff's Laws, this leads to the following \lhblue{equations for} a DGU bus $i \in \mathcal{D}$.
\begin{subequations} \label{eq_dc_dynamics}
\begin{gather}
l_i \dot{i}_i = -r_i i_i - v_i + u_i \\
c_i \dot{v}_i = i_i - i_i^o - g_iv_i
\end{gather}
\end{subequations}
where $l_i$, $r_i$ and $c_i$ are the filter inductance, resistance and capacitance, $g_i$ is a parasitic conductance, $i_i$ is the filter current, $v_i$ is the voltage at the output of the filter, $u_i$ is the voltage output of the buck converter (the control variable) and $i_i^o$ is the current injection from the neighbouring lines given by an expression analogous to \eqref{eq_kcl}.

In order to track a voltage setpoint (which can be broadcast by a centralised microgrid controller at semi-regular intervals to account for load changes and optimal power sharing within the network), each converter also contains an integrator of the form \cite{sadabadi_PlugandPlayRobustVoltage2018}
\begin{equation} \label{eq_dc_integrator}
\dot{\zeta}_i = v_i - v_i^{\mathrm{set}} + z_i i_i^o
\end{equation}
where $\zeta_i$ is the integrator state, $v_i^{\mathrm{set}}$ is the setpoint for converter $i \in \mathcal{D}$ and $z_i$ is a virtual impedance.

The dynamics \eqref{eq_dc_dynamics} and \eqref{eq_dc_integrator} together form the subsystem $i \in \mathcal{D}$, where the input from the line subsystem is $-i_i^o$ and the output is $v_i$, which can be written in the form \eqref{eq_linearised_sys}.

\begin{remark}
\lhblue{A virtual impedance can be used to enhance microgrid stability in the presence of destabilising loads \cite{lu_StabilityEnhancementBased2015}, and was also required for a passive controller of the form \eqref{eq_uistar} to be feasible in the case study considered, thus allowing for plug-and-play operation.} Using a virtual impedance will have an impact on the equilibrium voltage depending on the value of $Z_i$ and the equilibrium value of $i_i^o$. This can be accounted for in the value used for $v_i^{\mathrm{set}}$, which can be calculated and broadcast to each converter by a centralised controller responsible for maintaining power sharing in the microgrid after load changes.
\end{remark}

Meanwhile, we assume constant impedance loads at every $l \in \mathcal{L}$ with the following dynamics:
\begin{equation}
c_l \dot{v}_l = - g_l v_l  - i_l^o
\end{equation}
where $c_l$ is the load capacitance, $g_l$ is the load conductance and the rest of the notation has \lhblue{analogous meanings to those} in \eqref{eq_dc_dynamics}.  We have that $-i_l^o$ is the \lhblue{input and $v_l$ is the output for the load bus}.
We note that each $l \in \mathcal{L}$ is strictly passive with storage function $V_l(v_l) = \frac{1}{2} c_l v_l^2$.

Finally each line $(i,j) \in \mathcal{E}$ has the following dynamics:
\begin{equation}
l_{ij} \dot{i}_{ij} = -r_{ij} i_{ij} + v_{ij}
\end{equation}
where $l_{ij}$ and $r_{ij}$ are the line inductance and resistance, $i_{ij}$ is the line current, and $v_{ij}$ is the voltage difference between buses $i$ and $j$, given by an expression analogous to \eqref{eq_kvl}. The input from the combined DGU/load subsystems to each line system is $v_{ij}$ and the output is $i_{ij}$. We note that \lhteal{the dynamics associated with each} $(i,j) \in \mathcal{E}$ is strictly passive with storage function $V_{ij}(i_{ij}) = \frac{1}{2} l_{ij} i_{ij}^2$.

\begin{figure}[t]
\centering

\begin{tikzpicture}

\node[draw, circle,
			minimum size=0.6cm]
			(neg) at (0,0) {};
			\node at (neg.center){ $-1$};
\node[right=of neg,
			xshift=-0.5cm]
			(node1)  {};
\node[draw, circle,
			right=of node1,
			yshift= 0.5cm,
			xshift=-0.5cm,
			minimum size=0.6cm]
			(BD) {};
			\node at (BD.center){ $\calb_{\mathcal{D}}$};
\node[draw, circle,
			right=of node1,
			yshift= -0.5cm,
			xshift=-0.5cm,
			minimum size=0.6cm]
			(BL) {};
			\node at (BL.center){ $\calb_{\mathcal{L}}$};
\node[draw, rectangle,
			right=of BD,
			minimum width=1cm, 	minimum height=0.7cm]
			(sys_inv) {$\Sigma_{\mathrm{conv}}$};
\node[draw, rectangle,
			right=of BL,
			minimum width=1cm, 	minimum height=0.7cm]
			(sys_load) {$\Sigma_{\mathrm{load}}$};			
\node[draw, rectangle,
			below=of sys_inv,
			yshift= -0.25cm,
			minimum width=1cm, 	minimum height=1cm]
			(sys_line) {$\Sigma_{\mathrm{line}}$};
\node[draw, circle,
			right=of sys_inv,
			minimum size=0.6cm]
			(BtD)   {};
			\node at (BtD.center){ $\calb^T_{\mathcal{D}}$};
\node[draw, circle,
			right=of sys_load,
			minimum size=0.6cm]
			(BtL)   {};
			\node at (BtL.center){ $\calb^T_{\mathcal{L}}$};
\node[right=of BtL,
			yshift=0.5cm,
			xshift=-0.5cm]
			(node2)  {};
\node[right=of node2,
			xshift=-0.5cm]
			(node3)  {};
			
\draw (neg.east) -- (node1.center);
\draw[-stealth] (node1.center) |- (BD.west);
\draw[-stealth] (node1.center) |- (BL.west);
\draw[-stealth] (BD.east) -- (sys_inv.west)
	node[midway,above]{$w_{\mathrm{conv}}$};
\draw[-stealth] (BL.east) -- (sys_load.west)
	node[midway,above]{$w_{\mathrm{load}}$};
\draw[-stealth] (sys_inv.east) -- (BtD.west)
	node[midway,above]{$y_{\mathrm{conv}}$};
\draw[-stealth] (sys_load.east) -- (BtL.west)
	node[midway,above]{$y_{\mathrm{load}}$};
\draw (BtD.east) -| (node2.center);
\draw (BtL.east) -| (node2.center);
\draw (node2.center) -- (node3.center);
\draw[-stealth] (node3.center) |- (sys_line.east)
	node[pos=0.75,above]{$w_{\mathrm{line}}$};
\draw[-stealth] (sys_line.west) -| (neg.south)
	node[pos=0.25,above]{$y_{\mathrm{line}}$};

\end{tikzpicture}

\caption{Negative feedback interconnection of DC microgrid system, where we have a DGU subsystem $\Sigma_{\mathrm{conv}}$ with state $x_{\mathrm{conv}}$, load subsystem $\Sigma_{\mathrm{load}}$ with state $x_{\mathrm{load}}$ and a line subsystem $\Sigma_{\mathrm{line}}$ with state $x_{\mathrm{line}}$.}
\label{fig_interconnection_dc}
\end{figure}
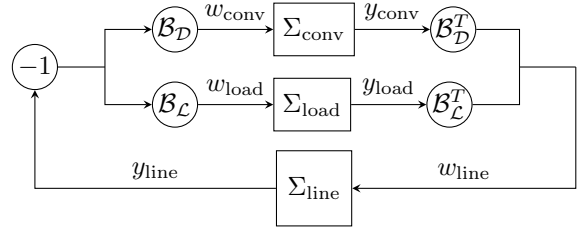

Now, to construct the closed-loop network system, the incidence matrix $\calb$ is divided into two sub-matrices: let $\calb_{\mathcal{D}}$ indicate the first $| \mathcal{D} |$ rows of $\calb$ and let $\calb_{\mathcal{L}}$ indicate the last $|\mathcal{L} |$ rows of the $\calb$. In addition, we can construct a composite converter system $\Sigma_{\mathrm{conv}}$ with state vector for the DGU \lhteal{subsystems} as $x_{\mathrm{conv}} = \left[ x_i \right]_{i \in \mathcal{D}}$, along with similar systems $\Sigma_{\mathrm{load}}$ and $\Sigma_{\mathrm{line}}$ with state vectors $x_{\mathrm{load}}$ and $x_{\mathrm{line}}$ for the load and line subsystems respectively. Let the notation $X_{\mathrm{conv}}$ indicate the matrix $X_{\mathrm{conv}} = \oplus_{i \in \mathcal{D}} X_i$, with similar definitions for $X_{\mathrm{load}}$ and $X_{\mathrm{line}}$. Then the full microgrid system is depicted in Figure \ref{fig_interconnection_dc} and can be described by $\dot{\hat{x}} = \hat{A}\hat{x} + \hat{B} + \hat{u}$ where $\hat{x} = \begin{bmatrix}
x_{\mathrm{conv}}^T & x_{\mathrm{load}}^T & x_{\mathrm{line}}^T
\end{bmatrix}^T$, $\hat{u} = \left[ u_i \right]_{i \in \mathcal{D}}$ and
\begin{subequations} \label{eq_microgrid_full}
\begin{gather}
\hat{A} =
\begin{bmatrix}
A_{\mathrm{conv}} & 0  &  - B_{\mathrm{conv}} \calb_{\mathcal{D}} C_{\mathrm{line}} \\
0 & A_{\mathrm{load}}  &  - B_{\mathrm{load}} \calb_{\mathcal{L}} C_{\mathrm{line}} \\
B_{\mathrm{line}} \calb_{\mathcal{D}}^T C_{\mathrm{conv}} & B_{\mathrm{line}} B_{\mathcal{L}}^T C_{\mathrm{load}}  &  A_{\mathrm{line}}
\end{bmatrix}  \\
\hat{B} = \lholive{
\begin{bmatrix}
B_{u \: \mathrm{conv}} &
0 &
0
\end{bmatrix}^T }
\end{gather}
\end{subequations}

Using Lemma \ref{lemma_lmi}, we can then construct matrices $P_i$ and $R_i$ for each $i \in \mathcal{D}$ such that each DGU subsystem is strictly passive from $-i_i^o$ to $v_i$ with storage function $V_i(x_i) = \frac{1}{2} x_i^T P_i x_i$ when the decentralised controller $u_i = -\frac{1}{2} R_i^{-1}B_{ui}^T P_i x_i$ is used to close the DGU subsystem control loop. This gives rise to a global storage function $\hat{V}(\hat{x}) = \frac{1}{2} \hat{x}^T \hat{P} \hat{x}$ with
\begin{equation} \label{eq_dc_P}
\hat{P} = \lholive{ \mathrm{blockdiag} \begin{bmatrix}
\oplus_{i \in \mathcal{D}} P_i &  \oplus_{l \in \mathcal{L}} c_l & \oplus_{(i,j) \in \mathcal{E}} l_{ij}
\end{bmatrix} }
\end{equation}
In addition, we have that the composite controller $\hat{u}$ is the optimal network-wide controller with respect to the LQR problem \eqref{opt_lqr} with cost matrices $R = \oplus_{i \in \mathcal{V}} R_i $ and (using $Q_{\calv}(R_i)$ as in \eqref{eq_q1q2q3})
\begin{equation}  \label{eq_dc_randq}
Q(R) =  \lholive{ \mathrm{blockdiag} \begin{bmatrix}
Q_{\calv}(R_i) & \oplus_{l \in \mathcal{L}} g_l &  \oplus_{(i,j) \in \mathcal{E}} r_{ij}
\end{bmatrix} }
\end{equation}

\subsubsection*{Simulation}
\begin{figure}[t]
\centering
\begin{tikzpicture}

\node[draw, circle,
			minimum size=0.1cm,
			align=center]
			(dgu1) at (0,0) {\tiny 1 DGU};
\node[draw, circle,
			right=of dgu1,
			minimum size=0.1cm,
			align=center]
			(dgu2) {\tiny 2 DGU};
\node[draw, circle,
			below=of dgu1,
			yshift=+0.75cm,
			minimum size=0.1cm,
			align=center]
			(load4) {\tiny 4 Load};
\node[draw, circle,
			below=of dgu2,
			yshift=+0.75cm,
			minimum size=0.1cm,
			align=center]
			(load5) {\tiny 5 Load};
\node[draw, circle,
			right=of load5,
			minimum size=0.1cm,
			align=center]
			(dgu3) {\tiny 3 DGU};

\draw[-stealth] (dgu1.south) -- (load4.north);
\draw[-stealth] (dgu1.east) -- (dgu2.west);
\draw[-stealth] (dgu2.south) -- (load5.north);
\draw[-stealth,dotted] (load5.east) -- (dgu3.west);

\end{tikzpicture}
\caption{Graph structure of the DC microgrid used in simulation. The line between DGU 3 and Load 5 is only connected at $t = 3$ seconds.}
\label{fig_microgrid}
\end{figure}
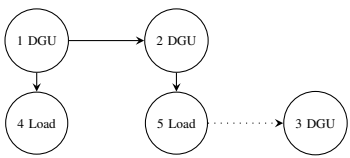

\begin{figure}[t]
\centering
\includegraphics[width=0.45\textwidth]{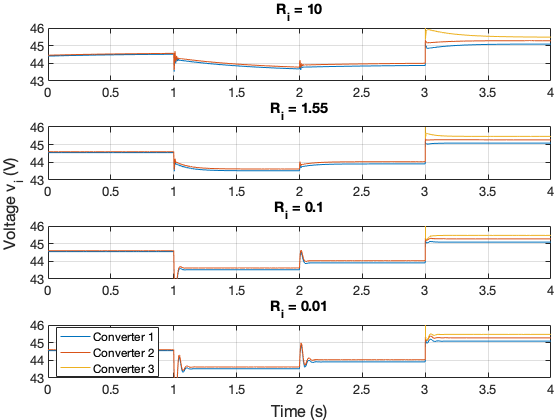}
\caption{The voltage $v^i$ generated by each converter in a \lhgreen{five}-bus system containing \lhgreen{three} DGU buses and two load buses for various values for $\bar{R}_i$.}
\label{fig_compare_R_dc}
\end{figure}

The above results have been tested numerically on a simple microgrid composed of three DGU buses and two load buses (see Figure \ref{fig_microgrid}). The three DGU buses have identical converters, with parameters given in Table \ref{tab_dc} (when no specific index is provided, the value applies across all relevant buses and lines). \lhorange{The LMIs \eqref{opt_lmi} were first solved using the parameters for each DGU bus, then the resulting controllers with various values for $\bar{R}_i$ were tested in simulations to evaluate performance against a series of disturbances. The results can be found in Figure \ref{fig_compare_R_dc}. Here, the values of $g_4$ was changed to \SI{0.15}{S} at $ = 1$ second, then the value of $g_5$ was changed to \SI{0.03}{S} at $t = 2$ seconds. Finally, the link between DGU 3 and Load 5 was connected at $t = 3$ seconds to test the plug-and-play property.}

\begin{table}[t]
\caption{DC Microgrid parameters (adapted from \cite{DecentralizedScalableApproachtucci2015})}
\begin{scriptsize}
\begin{center}
\begin{tabular}{|cc|cc|}
\hline
\multicolumn{4}{|c|}{\textbf{Converter Parameters}} \\
\hline
$r_i$ &  \SI{0.2}{\ohm} & $v_1^{\mathrm{set}}$ & \SI{48}{V}   \\
$c_i$ &  \SI{2.2}{mF} & $v_2^{\mathrm{set}}$ & \SI{47.8}{V}  \\
$l_i$ &  \SI{1.8}{mH} & $v_3^{\mathrm{set}}$  & \SI{48.1}{V}  \\
$g_i$ &  $\frac{1}{100}$\SI{}{S} & $z_i$  & \SI{1}{\ohm}    \\
\hline
\multicolumn{4}{|c|}{\textbf{Load Parameters}} \\
\hline
$g_4$ &  \SI{0.1}{S} & $g_5$ &  \SI{0.05}{S}  \\
$c_4$ & \SI{70}{\micro F} & $c_5$ &  \SI{70}{\micro F}  \\
\hline
\multicolumn{4}{|c|}{\textbf{Line Parameters}} \\
\hline
$r_{ij}$ & \SI{0.05}{\ohm} & $l_{ij}$ & \SI{2.1}{\micro H} \\
\hline
\end{tabular}
\label{tab_dc}
\end{center}
\end{scriptsize}
\end{table}

\lhorange{The values of $R_i$, $Q(R)$ and $P_i$ were found when \eqref{opt_lmi} was}  ran as an optimisation problem to find the maximal $R_i$, as outlined in Remark \ref{rem_opt_problem}. In addition, the constraint \eqref{eq_lambda} was included in the optimisation problem with $\lambda_{\lhmag{i}} = -8$ to find solutions that yielded a reasonable settling time. The value of $R_i$ obtained was 1.55, from which the corresponding values of the global matrices $R$, $\hat{P}$ and $Q(R)$ could be calculated using \eqref{eq_dc_P} and \eqref{eq_dc_randq}.

Following this, tuning could be achieved by altering the value of $\bar{R}_i$, with values less than 1.55 simultaneously generating a strictly passive bus and an optimal controller. As ${R}_i = 1.55$ is the maximal value that can satisfy all conditions in Theorem \ref{theorem_combined} simultaneously, $\bar{R}_i > 1.55$ \lhblue{cannot both passivate the bus dynamics and have an inverse optimal interpretation} (but may still lead to stable dynamics).

As can be seen \lhorange{in Figure \ref{fig_compare_R_dc}}, the controller quickly returns the microgrid voltage to equilibrium after each disturbance. \lhblue{In addition, the change in microgrid topology at $t=3$ seconds \lhorange{is} easily handled}. New equilibria are established after each load change due to the effect of the virtual impedance $z_i$, however deviations are relatively small and can be corrected by a new voltage setpoint for each converter.  As expected, lower values of $\bar{R}_i$ allow equilibrium to be established more quickly as the cost of control is lowered and the corresponding cost of non-equilibrium values for the states is increased in the matrix $Q(\bar{R}_i)$. However, lower values of $\bar{R}_i$ result in a controller with higher gains, which can have the effect of an overly aggressive action, which can cause large voltage overshoots as in the case of $\bar{R}_i = 0.1$ and $\bar{R}_i = 0.01$ in Figure \ref{fig_compare_R_dc}. Nevertheless, it is clear that the framework \lhblue{described} here offers a practical approach towards designing a controller for a DC microgrid system which allows for plug-and-play \lhblue{capabilities} with guaranteed stability, and also allows for performance to be quantified via an optimal \lhblue{control problem for the network}.

\section{Conclusion} \label{sec_conclusion}
In this paper, we have introduced a set of sufficient conditions whereby a decentralised passivity-based controller can be measured with respect to a global network performance metric. A method to design such a decentralised control system was provided, which has the dual benefit of simultaneously rendering the closed-loop system \lhred{at each node} strictly passive, thus 
\icl{ensuring stability and facilitating a} plug-and-play \lhblue{\icl{operation}}, while also ensuring that the synthesised controller is the solution to a network-wide optimal control problem. This was achieved by combining passivity theory with inverse optimal control theory to produce a set of constraints on the local node dynamics in the form of linear matrix inequalities that depend only on local node parameters. The approach was then verified using a DC microgrid case study.

\appendix
\subsection{Proof of Theorem \ref{theorem_combined}} \label{app_proof_t2}
The proof relies on direct application of Theorem \ref{theorem_inverse_opt} to the network system \eqref{eq_full_network}. Firstly take $\hat{V}(\hat{x}) = V_{\calv}(x_{\calv}) + V_{\cale}(x_{\cale})$ where $V_{\calv}(x_{\calv}) = \sum_{i \in \mathcal{V}} V_i(x_i)$ and $V_{\cale}(x_{\cale}) = \sum_{j \in \mathcal{E}} V_j(x_j)$ so that \lhblue{$
\nabla \hat{V}(\hat{x}) = \begin{bmatrix}
\nabla^T V_{\calv}(x_{\calv}) & \nabla^T V_{\cale}(x_{\cale})
\end{bmatrix}^T$
where $\nabla V_{\calv}(x_{\calv}) = \left[ \nabla V_i(x_i) \right]_{i \in \calv}$ and $ \nabla V_{\cale}(x_{\cale}) = \left[ \nabla V_j(x_j) \right]_{j \in \cale}$}.

Now, the optimal network wide controller is of the form $\hat{u}^{\star}(\hat{x}) = -\frac{1}{2} R^{-1}\hat{B}^T \nabla \hat{V}$, \lhteal{with $\hat{B}$ given in \eqref{eq_network_linear_full}}. As $\hat{B}^T$ has a block-diagonal structure, it is clear that in order for $\hat{u}^{\star}(\hat{x})$ to be decentralised (i.e. for each $u_i, i \in \mathcal{V}$ to depend only on local $x_i$), then $R$ must also have a block-diagonal structure. If $R$ is chosen as $R = \oplus_{i \in \mathcal{V}} R_i$, then explicit calculation of $\hat{u}^{\star}$ gives $\hat{u}^{\star}(\hat{x}) = \left[ - \tfrac{1}{2} R_i^{-1}B_{ui}^T(x_i) \nabla V_i(x_i) \right]_{i \in \mathcal{V}} = \left[ u_i \right]_{i \in \mathcal{V}}$,
proving that $u_i$ is indeed composed of the individual decentralised controllers.

In addition, \lhteal{
\icl{\eqref{eq_state_cost_general} 
gives}}
\begin{multline} \label{eq_q_expansion_1}
\nabla^T \hat{V}(\hat{x}) \hat{f}(\hat{x}) =  \nabla^TV_{\calv}(x_{\calv}) f_{\calv}(x_{\calv}) + \nabla^T V_{\cale}(x_{\cale}) f_{\cale}(x_{\cale}) \\ - \nabla^T V_{\calv}(x_{\calv}) B_{\calv}(x_{\calv}) \calb_p y_{\cale} + \nabla^T V_{\cale} B_{\cale}(x_{\cale})\calb_p^Ty_{\calv}
\end{multline}
Now, because each $i \in \mathcal{V}$ and $j \in \mathcal{E}$ is \lhviolet{strictly} passive, \icl{the storage functions can be chosen such that} $\nabla^TV_k(x_k) B_k(x_k) = h_k^T(x_k), \forall k \in \mathcal{V} \cup \mathcal{E}
$
\lhteal{\cite{hill_StabilityNonlinearDissipative1976}} which implies
$
\nabla^T V_{\calv}(x_{\calv})B_{\calv}(x_{\calv}) = y_{\calv}^T $ and $
\nabla^T V_{\cale}(x_{\cale})B_{\cale}(x_{\cale}) = y_{\cale}^T
$
Therefore, the last two terms in \eqref{eq_q_expansion_1} cancel. In addition, \lhblue{
$
\nabla^T \hat{V}(\hat{x})\hat{B}(\hat{x})\hat{u} $ $=$ $ - \sum_{i \in \mathcal{V}} \tfrac{1}{2} \nabla^T V_i(x_i) B_{ui}(x_i)  R_i^{-1}B_{ui}^T(x_i)\nabla V_i(x_i)
$ and
$
\lVert \hat{u} \rVert^2_R$ $=$ $\sum_{i \in \mathcal{V}}  \tfrac{1}{4}\nabla^TV_i(x_i) B_{ui}(x_i) R_i^{-1}B_{ui}^T(x_i)\nabla V_i(x_i)
$.}
Combining the above gives $q(\hat{x},R)$ as in \eqref{eq_state_cost_network}.

It remains to show that condition \eqref{eq_condition} is satisfied, which is equivalent to $q(\hat{x},R) > 0$. As we have \eqref{eq_condition_local} by assumption for all $i \in \mathcal{V}$, then we have that the summation across terms involving $\mathcal{V}$ in \eqref{eq_state_cost_network} is positive. For each $j \in \mathcal{E}$, we have $\nabla^TV_j(x_j)f_j(x_j) \leq 0$ by \lhviolet{strict} passivity, so the summation of terms involving $\mathcal{E}$ in \eqref{eq_state_cost_network} is non-negative, \lhteal{as required}.

Finally, \lhteal{asymptotic} stability \lhblue{follows directly from the fact that the network is composed of the negative feedback interconnection of \lhviolet{two strictly passive subsystems} (see Section \ref{sec_network_model}), leading to the Lyapunov function $\hat{V}(\hat{x}) = \sum_{i \in \calv} V_i(x_i) + \sum_{j \in \cale} V_j(x_j)$ \lhviolet{ \cite[Theorem 6.3]{khalil_nonlinear_2013}}}.

\subsection{Proof of Lemma \ref{lemma_lmi}} \label{app_proof_l3}
Firstly, let the local state feedback controller matrix be of the form $u_i = K_ix_i$ with $K_i = -\frac{1}{2}R_i^{-1} B_{ui}^T P_i$, where $i \in \mathcal{D}$, and let $S_i = -\frac{1}{2}R_i^{-1}$ so that $K_i = S_i B_{ui}^T P_i$. \lhgreen{Then}:
\begin{enumerate}
\item Substitution of $K_i$ above into \eqref{eq_solve_LMI} gives
$
\begin{bsmallmatrix}
\mathcal{A} + \Gamma_i^{-1} & \mathcal{C} \\
\mathcal{C}^T & 0
\end{bsmallmatrix} \quad \leq 0
$
where $\mathcal{A} = A_i^TP_i + P_iA_i + 2P_iB_{ui}S_iB_{ui}^TP_i$ and $\mathcal{C} = P_i B_{i}- C_i^T$. Letting $Y_i = P_i^{-1}$ and pre- and post-multiplying \eqref{eq_solve_LMI} by $\begin{bsmallmatrix} Y_i & 0 \\ 0 & I \end{bsmallmatrix}$ and utilisation of the Schur \lhteal{complement} yields LMI \eqref{const_passive}. Therefore, if \eqref{const_passive} is satisfied, we have that $i \in \mathcal{V}$ is strictly passive (\eqref{eq_pass_condition_local} is satisfied) and the controller is of the form \eqref{eq_control_law_local} with $V_i(x_i) = \frac{1}{2}x_i^T P_i x_i$.

\item The local condition \eqref{eq_condition_local} for the linear system reduces to constraint \eqref{const_q_pos_def}, so using the solutions to \eqref{opt_lmi} in \eqref{eq_uistar} satisfies the conditions of Theorem \ref{theorem_combined}. Following this, direct calculation of $\hat{u}^{\star} = \left[ u_i \right]_{i \in \mathcal{V}}$, $R = \oplus_{i \in \calv} R_i$ and $Q(R)$ (as per \eqref{eq_state_cost_network}) yields the optimal controller \lhteal{\eqref{eq_ustar}} and associated cost matrices given by \lhteal{\eqref{eq_RandQ}}.

\item \lhteal{Asymptotic} stability also follows from Theorem \ref{theorem_combined}.
\end{enumerate}

\subsection{Proof of Lemma \ref{lemma_tuning}} \label{app_proof_l4}
Take $S_i = - \frac{1}{2}R_i^{-1}$ and $\bar{S}_i = - \frac{1}{2}\bar{R}_i^{-1}$. We have $\bar{R}_i \leq R_i$, which implies $\bar{S}_i \leq S_i$, so let $\bar{S}_i = S_i - \Delta$, where $\Delta \geq 0$.
\begin{enumerate}
\item To check for passivity, we can investigate if \eqref{eq_solve_LMI} is satisfied with $\bar{S}_i$ and $P_i$. We have
\begin{equation*} \lhgreen{
\begin{aligned}
\bar{\mathcal{A}} + \Gamma_i^{-1} & =  A_i^TP_i + P_iA_i + 2P_iB_{ui}S_iB_{ui}^TP_i + \Gamma_i^{-1} \\ & \qquad \lhteal{-} 2P_iB_{ui} \Delta B_{ui}^TP_i \quad \lhteal{\leq} \quad 0
\end{aligned} }
\end{equation*}
where the first set of terms is \lhteal{negative} semi-definite by virtue of the fact that $R_i$ and $P_i$ satisfy \eqref{const_passive}. Also \lhteal{the blocks in \eqref{eq_solve_LMI} containing $P_iB_i - C_i^T$ remain} unchanged, therefore \eqref{eq_solve_LMI} is satisfied with $\bar{S}_i$ and $P_i$ and the local subsystem is passive.

\item We have that $Q_{\calv}(\bar{R}_i) \geq Q_{\calv}(R_i) > 0$ (as $\bar{R}_i^{-1} \geq R_i^{-1}$) and $Q_{\cale}$ remains unchanged, so $Q(\bar{R}) \geq 0$. Therefore, Theorem \ref{theorem_inverse_opt} can be applied, giving $\bar{u}$ as the optimal controller to \eqref{opt_lqr} with cost matrices $\bar{R}$ and $Q(\bar{R})$.
\end{enumerate}

\bibliographystyle{IEEEtran}
\bibliography{bibliography}
\end{document}